\documentclass{svproc}

\usepackage{url}
\usepackage{graphicx}

\usepackage{amsmath,amssymb,xcolor}

\def\d{\mathrm{d}}
\def\eqref#1{(\ref{#1})}
\def\RR{\mathbb{R}}

\newtheorem{assumption}{Assumption}

\begin{document}
\mainmatter              %
\title{On implicit interpolation models for nonlinear anisotropic magnetic material behavior}
\titlerunning{Implicit models for magnetic anisotropy}  %
\author{Herbert Egger~\inst{1,2} \and Michael Mandlmayr\inst{2}
}
\authorrunning{Herbert Egger et al.} %
\tocauthor{Herbert Egger and Michael Mandlmayr}
\institute{Johannes Kepler University Linz, Austria
\and
Johann Radon Institute for Computational and Applied Mathematics, Linz, Austria}

\maketitle              %

\begin{abstract}
Implicit models for magnetic coenergy have been proposed by Pera et al. to describe the anisotropic nonlinear material behavior of electrical steel sheets. This approach aims at predicting magnetic response for any direction of excitation by interpolating measured of B--H curves in the rolling and transverse directions. In an analogous manner, an implicit model for magnetic energy is proposed. We highlight some mathematical properties of these implicit models and discuss their numerical realization, outline the computation of magnetic material laws via implicit differentiation, and discuss the potential use for finite element analysis in the context of nonlinear magnetostatics.
\keywords{anisotropic materials, electrical steel sheets, nonlinear magnetostatics, finite element analysis}
\end{abstract}
\section{Introduction}

The efficient operation of power transformers and electric motors critically relies on the excellent magnetic permeability and low electrical conductivity of modern electrical steel sheets. 
Accurate modeling of the nonlinear magnetic behavior of these materials therefore is a key ingredient for the systematic analysis and optimization of such devices. 
Besides magnetic saturation also the anisotropic behavior of grain oriented and non-oriented electrical steel sheets has been recognized as an important aspect for modeling and simulation. Various approaches to describe anisotropic permeability tensors and their use in finite element analysis have been investigated in the literature; see~\cite{ferrara21,higuchi14,jiang17,martin16} for some recent examples and further references.  
\paragraph{General approach.}
Following Silvester and Gupta \cite{silvester91}, the anhysteretic anisotro\-pic nonlinear response of a soft magnetic material can be completely described by an energy functional $w=w(b)$ or a respective coenergy functional $w_*=w_*(h)$. 
Equivalence of the two formulations can be established by convex duality \cite{rockafellar70}.
The relation between magnetic field $h$ and magnetic flux $b$ is then given by 
\begin{equation} \label{eq:1}
h(b)=\nabla_b w(b) \qquad \mbox{resp.} \qquad 
b(h) = \nabla_h w_*(h),
\end{equation}
and the Jacobians of these vector valued relations 
\begin{equation} \label{eq:2}
\nu'(b) =\frac{\d h(b)}{\d b} \qquad \mbox{resp.} \qquad \mu'(h) = \frac{\d b(h)}{\d h}
\end{equation}
define the differential reluctivity and permeability tensors, which are inverse to each other.
Since they correspond to the Hessian matrices for the underlying energy and coenergy functionals, respectively, one can see that they are symmetric and positive definite, if the underlying functionals are strongly convex. 

The description of appropriate energy or coenergy functionals thus allows to completely describe the anhysteretic response of an anisotropic nonlinear material, providing all the necessary information required for nonlinear finite element analysis using vector or scalar potential formulations \cite{kameari06,meunier08,vandesande04}.
As mentioned in \cite{silvester91}, the true energy or coenergy functional can be approximated by a multivariate spline. While this provides a flexible approximation and a convenient implementation of the material laws \eqref{eq:1}--\eqref{eq:2}, the fitting of a multivariate spline requires measurements of B--H curves for many directions, which may be hardly available in practice.
To reduce the number of degrees of freedom in the fitting process, the use of a certain parametric form of the energy functional have been proposed in \cite{martin16}. 
Alternative parametrizations have been used previously for the modeling of anisotropic permeabilities or reluctivities \cite{dinapoli83,kameari06}.

\paragraph{Implicit interpolation.}
In \cite{pera93}, Pera et al. proposed a procedure to predict anisotropic response in arbitrary directions from measurements of B--H curves in the rolling and transverse direction only. The approach relies on characterizing contours of equal coenergy $w_*$ via a nonlinear equation of the form
\begin{equation} \label{eq:3}
\Big|\frac{h_1}{\hat h_1(w_*)}\Big|^{n(w_*)} + \Big|\frac{h_2}{\hat h_2(w_*)}\Big|^{n(w_*)} = 1. 
\end{equation}
The two functions $\hat h_1(w_*)$, $\hat h_2(w_*)$ result from inverting the scalar functional relations $w_*=w_*(h_1,0)$, $w_*=w_*(0,h_2)$ for the coenergy along the rolling and transverse direction, which can be obtained from B--H measurements for these orientations \cite{fryskowksi08,pera94}, and the exponent $n(w_*)$ influences the shape of the contours of equal coenergy \cite{pera93}.
For $n(w_*)=2$, these correspond to ellipses \cite{dinapoli83}, but $n(w_*)>2$ is required to reproduce the typical anisotropic behavior \cite{pera94}.
The equation \eqref{eq:3} implicitly defines the coenergy $w_*=w_*(h_1,h_2)$ as a function of the two field components by some sort of implicit interpolation between the corresponding values in rolling and transverse direction. 
Numerical differentiation was used in~\cite{pera93,pera94} to compute approximations of \eqref{eq:1} and \eqref{eq:2}, which can be used in finite element simulations based on a magnetic scalar potential. 

In a completely analogous manner, we can define a relation of the form 
\begin{equation} \label{eq:4}
\Big|\frac{b_1}{\hat b_1(w)}\Big|^{p(w)} + \Big|\frac{b_2}{\hat b_2(w)}\Big|^{p(w)} = 1 
\end{equation}
to describe contours of equal energy and thus implicitly define an energy functional $w=w(b_1,b_2)$, which is the form suitable for simulations based on a magnetic vector potential~\cite{meunier08}. Here the exponent $p(w)<2$ is needed to obtain realistic behavior at intermediate directions.
Although obvious, this approach seems not to have been discussed extensively in the literature so far. 
A somewhat similar idea was proposed in \cite{biro10} and used in \cite{chwastek13}. These works however consider $\hat b_1(\cdot)$, $\hat b_2(\cdot)$ and $p(\cdot)$ as functions of $|h|$ instead of $w$. In that case, an additional relation is required to obtain $h=h(b)$, and it is unclear, if there exists a corresponding energy functional $w(b)$ such that $h(b)=\nabla_b w(b)$. Consistency of this material model with basic thermodynamic relations is thus unclear.  

\paragraph{Outline and main contributions.}
In this paper, we study the two implicit aniso\-tropic material models \eqref{eq:3} and \eqref{eq:4} from a mathematical point of view. 
In particular, we study the well-posedness and the efficient evaluation of the energy and coenergy functionals $w(b)$ and $w_*(h)$, and we indicate how to practically realize \eqref{eq:1} and \eqref{eq:2} without using numerical differentiation. 
In addition, we study additional mathematical properties that are required to guarantee well-posedness of the corresponding finite element models as well as convergence of iterative solvers for the nonlinear systems, namely, \emph{smoothness}, \emph{convexity}, and \emph{coercivity} of the energy resp. coenergy functionals; see \cite{engertsberger23,heise94} for results in this direction.
Numerical examples are presented for illustration and some conclusions of our investigations are summarized at the end of the manuscript.

\section{Mathematical properties}

Following \cite{biro10,pera93}, we consider a planar two-dimensional setting. 
The coordinate system is oriented such that the two axis are aligned with the rolling and transverse direction. 
Throughout the following discussion, we focus on the coenergy setting \eqref{eq:3}, but with similar arguments, the energy model \eqref{eq:4} can be analyzed as well. 
Comments on this generalization are given at the end of the section.  

\subsection{Main assumptions}

For our analysis, we use some elementary properties of the functions $\hat h_i(\cdot)$ appearing in \eqref{eq:3}, which can be deduced from the following basic assumptions on B--H curves in the principal directions available from measurements~\cite{fryskowksi08,pera94}. 
\begin{assumption} \label{ass:1}
Along the principal axes, the magnetic flux $b$ and the magnetic field $h$ are colinear, i.e., 
\begin{equation}
b(h_1,0)=\binom{b_1(h_1)}{0} \qquad \text{and} \qquad b(0,h_2)=\binom{0}{b_2(h_2)},
\end{equation}
with scalar functions $b_i(h_i)$, $i=1,2$. 
These functions are smooth, uniformly monotone, globally Lipschitz, and antisymmetric, i.e., 
\begin{equation}
0 < \underline \beta \le b_i'(h_i) \le \overline \beta < \infty 
\qquad \mbox{and} \qquad 
b_i(-h_i)=-b_i(h_i).
\end{equation}
\end{assumption}
These assumptions match with the observations made in measurements and with the usual conditions used for the analysis of isotropic materials; see e.g. \cite{dinapoli83,pera94}. 

\subsection{Preliminary considerations}

By integration, we may define $w_{*,i}(h_i) = \int_0^{h_i} b_i(s) \, ds$ which are strictly monotone for positive arguments, since the integrands are positive. 
We may thus define  corresponding monotonically increasing inverse funtions $\hat h_i(w_*)$, which are monotonically and convex, and will serve as building blocks for the model \eqref{eq:3}. 

\subsection{Well-posedness and basic properties}

As a next step, we show that under the previous assumptions, the model \eqref{eq:3} in fact defines a proper coenergy functional $w_*(h_1,h_2)$ with some useful properties. 
For the following statement, we require the exponent $n(w_*)$ to be constant, i.e., independent of $w_*$. We will comment on the general case in remarks later on. 
\begin{lemma} \label{lem:2}
Let Assumption~\ref{ass:1} hold and $n(w_*)=n > 0$ constant. 
Then for any $(h_1,h_2) \neq 0$, equation \eqref{eq:3} has a unique solution $w_*=w_*(h_1,h_2)>0$.
Moreover, $w_* \to 0$ mit $(h_1,h_2)\to 0$ and $w_* \to \infty$ with $|(h_1,h_2)| \to \infty$.
The function $w_*:\RR^2 \setminus 0 \to \RR_+$, $(h_1,h_2) \to w_*(h_1,h_2)$ is continuous, coercive, and symmetric in both arguments, and can be extended continuously to $\RR^2$ by $w_*(0,0)=0$.
\end{lemma}
The proof relies on the monotonicity of the right hand side in equation~\eqref{eq:3} with respect to $w_*$, which can be guaranteed under the assumptions of the lemma.
For variable $n(w_*)$, monotonicity, and hence uniqueness of a solution $w_*$ to \eqref{eq:3} does in general not hold, which can be shown by examples.
\begin{remark}
The evaluation of $w_*=w_*(h_1,h_2)$ requires the solution of the scalar nonlinear equation \eqref{eq:3}, which can be accomplished by standard methods, e.g., bisection, regula-falsi, Newton. 
Under the assumptions of Lemma~\ref{lem:2}, convergence of these algorithms can be guaranteed and thus the evaluation of $w_*(h_1,h_2)$ can be realized efficiently and reliably. 
\end{remark}

\subsection{Differentiability}

In order to utilize \eqref{eq:1} and \eqref{eq:2} for finite element analysis, we need to evaluate these equations and hence have to ensure sufficient smoothness of the coenergy function $w_*(h_1,h_2)$, which is only defined implicitly via equation~\eqref{eq:3}.  
According to the previous remark, we consider only the choice $n(w_*)=n$ constant.
\begin{lemma} \label{lem:3}
Let Assumption~\ref{ass:1} hold and $n(w_*)=n \ge 1$ constant. 
Then the function $w_*(h_1,h_2)$ defined by \eqref{eq:3} is continuously differentiable. 
For exponent $n \ge 2$, the function $w_*(h_1,h_2)$ is twice continuously differentiable on $\RR^2 \setminus (0,0)$.
\end{lemma}
\begin{proof}
By inserting $w_*=w_*(h_1,h_2)$ in \eqref{eq:3} and formal differentiation with respect to $h_i$, one otains
\begin{align*}
D(h_1,h_2,w_*) \, \partial_i w_* = \left(\frac{h_i^2}{\hat h_i^2(w_*)}\right)^{\frac{n}{2}-1} \frac{2 h_i}{\hat h_i^2(w_*)}.
\end{align*}
Under the assumptions of the lemma, the  discriminant $D(h_1,h_2,w_*)$ can be seen to be strictly positive. This allows to compute the derivatives $\partial_i w_*(h_1,h_2)$ of the coenergy, once the solution $w_*$ of \eqref{eq:3} has been determined. 
In a similar manner, one can derive a formula for the second derivatives. 
\qed
\end{proof}

\begin{remark} 
The proof yields formulas for computing the derivatives, i.e., for evaluating the material laws \eqref{eq:1}--\eqref{eq:2}, which are required for the implementation in finite element analysis. 
In particular, no numerical differentiation as in \cite{pera93,pera94} is required. 
Similar formulas can also be derived for variable exponent $n(w_*)$, although differentiability of the discriminant cannot be ensured a-priori. 
\end{remark}

\begin{remark}
An exponent $n>2$ is required to obtain a physically reasonable behavior, i.e., a direction of hard magnetization which is not one of the principal axes \cite{pera93,pera94}.
In this case, twice continuous differentiability of the coenergy function $w_*(h_1,h_2)$ can be expected. 
For $n < 2$, the second derivatives are well-defined for $h_1, h_2 \ne 0$, which however may blow up close to the axes. %
\end{remark}

\subsection{Convexity}

Besides smoothness, the convexity is another important property of the coenergy functional $w_*(h_1,h_2)$ which is required to ensure, e.g., positive definitness of the incremental permeability $\mu'(h) = \frac{\d b(h)}{\d h}$, thermodynamic consistency of the respective material model, well-posedness of corresponding finite element analyses, and convergence of iterative solution methods \cite{engertsberger23,krause12}.
To begin with, let us briefly summarize some important properties of convex functions. 

\begin{remark}
Let $w_* : \RR^2 \to \RR$ be convex. Then \\[1ex]
(i) The levelsets $L(s) := \{(h_1,h_2) : w_*(h_1,h_2)=s\}$ are convex. \\[1ex]
(ii) The scalar functions $w_*(h_1,0)$ and $w_*(0,h_2)$ are convex. \\[1ex]
(iii) If $w_*$ is strictly convex and twice continuously differentiable at $(h_1,h_2)$, then the Hessian $H w_*(h_1,h_2)$ is symmetric positive definite.  \\[1ex]
(iv) If $w_*$ is twice continuously differentiable and its Hessian symmetric and positive definite everywhere, then $w_*$ is strictly convex.
\end{remark}

From the particular construction of the implicit coenergy model \eqref{eq:3}, one can immediately deduce the following basic properties. 

\begin{lemma}
Let Assumption~\ref{ass:1} hold and $1 \le n(w_*) \le \infty$. Then the implicit coenergy model \eqref{eq:3} satisfies the conditions (i)--(ii)  of the previous remark. 
\end{lemma}
Unfortunately, the remaining properties (iii)--(iv) are not so obvious. 
For the case of variable exponent $n(w_*)$, examples which lack global convexity can again be constructed explicitly, and even in the case of constant exponent $n(w_*)=n$ it seems highly non-trivial to guarantee convexity a-priori. 
We however can identify a special case, for which convexity can still be proven rigorously. 
\begin{lemma} \label{lem:5}
Let Assumption~\ref{ass:1} hold and $n(w_*)=n \ge 1$ be constant. Further assume that $\hat h_1(w_*) = \lambda \hat h_2(w_*)$ for some $\lambda > 0$. Then the implicit coenergy function $w_*(h_1,h_2)$ defined by \eqref{eq:3} is convex.  
\end{lemma}
\begin{proof}
Inserting the relation $\hat h_1(w_*)=\lambda \hat h_2(w_*)$ into \eqref{eq:3} leads to
\begin{align*}
1 
=\left(\lambda^n |h_1|^n + |h_2|^n\right) / \hat h_2^2(w_*)^{n/2}.
\end{align*}
Using Assumption~\ref{ass:1} one can see that $\psi(w_*) :=\hat h_2^2(w_*)^{n/2}$ is strongly monotone and invertible on $\RR^+$. We can thus determine  $w_*$ explicitly by 
\begin{align} \label{eq:special}
w_*(h_1,h_2)=\psi^{-1}\left(\lambda^n |h_1|^n + |h_2|^n\right),
\end{align}
and this function is convex by construction and the previous arguments. 
\qed
\end{proof}

\begin{remark}
The special form \eqref{eq:special} of the coenergy functional yields a suitable generalization of the models $w_*(h_1,h_2)=f(|h_1|^2+|h_2|^2)$ usually employed in the isotropic case. Even the analysis for a variable exponent $n(w_*)$ seems feasible, since the model \eqref{eq:special} is explicit.
\end{remark}

\subsection{The corresponding energy model} \label{wec:2.6}
All assertions obtained for the coenergy model \eqref{eq:3} immediately translate also to the energy model \eqref{eq:4}. 
There, however, an exponent $p(w)<2$ is required to obtain the physically reasonable behavior; see \cite{biro10,martin16}. 
The second derivatives of the resulting energy functional $w(b_1,b_2)$ will have singularities at the principal axis, which should be taken into account carefully, if this model is used in finite element analysis.
This may also explain some of the numerical difficulties observed in \cite{biro10,chwastek13} for the implementation of a related model.

\section{Numerical validation}

In the following, we briefly discuss some typical examples that highlight the main features of the anisotropic implicit (co)energy models discussed above.

\paragraph*{Example~1.}
We assume that the magnetic response in the two principal directions is linear, i.e.,  $\hat h_1(w_*)=c_1\sqrt{2 w^*}$ and $\hat h_2(w_*)=c_2\sqrt{2 w^*}$.
From the calculations used in the proof of Lemma~\ref{lem:5}, we see that the coenergy defined implicitly via \eqref{eq:3} with constant exponent $n$ is then given by
\begin{align} \label{eq:ex1}
w^*(h_1,h_2)=\tfrac{1}{2}\|(h_1/c_1,h_2/ c_2)\|_n^2.
\end{align}
Here $\|(a,b)\|_p=(|a|^p+|b|^p)^{1/p}$ denotes the $\ell^p$-norm of the vector $(a,b)$.
For the following examples, we set $c_1=2$ and $c_2=1$, which thus correspond to the transverse and rolling direction, respectively.

\begin{figure}
    \centering
    \includegraphics[width=0.45\textwidth]{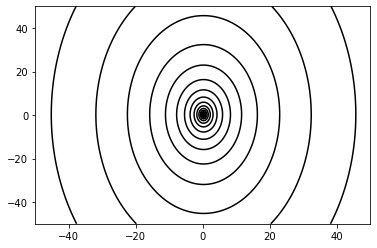}
    \includegraphics[width=0.45\textwidth]{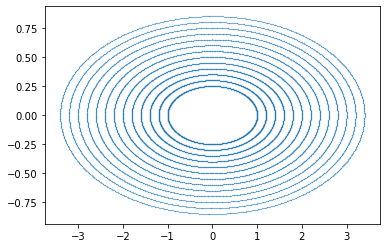}
    \caption{Contours of equal coenergy for the model  \eqref{eq:ex1} with $n=2$ (left), and locus of the $(h_1,h_2)$-curves (right) for constant induction $|(b_1,b_2)|$ computed by relation \eqref{eq:1}.}
    \label{fig:lin}
\end{figure}

\paragraph{Linear anisotropic model.}
For exponent $n=2$, the coenergy is a quadratic functions of $(h_1,h_2)$, i.e., the model is linear but still anisotropic , since $c_1 \ne c_2$.
In Figure~\ref{fig:lin}, we depict the contours of equal coenergy and the locus of the $(h_1,h_2)$-curves, that correspond to circles of constant strength of induction $|(b_1,b_2)|$.
In this example, the direction of hard magnetization is the transverse direction. 

\begin{figure}
    \centering
    \includegraphics[width=0.45\textwidth]{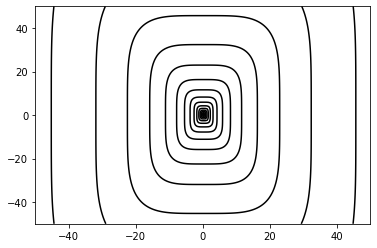}
    \includegraphics[width=0.45\textwidth]{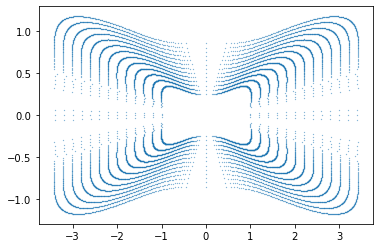}
    \caption{Contours of equal coenergy for the model  \eqref{eq:ex1} with $n=13/3$ (left), and locus of the $(h_1,h_2)$-curves (right) for constant induction $|(b_1,b_2)|$.}
    \label{fig:coenergy}
\end{figure}

\paragraph{Coenergy setting.}
For a second test, we choose the exponent $n=13/3>2$ in equation \eqref{eq:ex1}, which mimicks the physically relevant setting in the coenergy framework~\eqref{eq:3}. 
In Figure~\ref{fig:coenergy}, we again plot the contours of equal coenergy and the locus of the $(h_1,h_2)$ curves for constant induction. 
As intended, the hard direction now is different from the principal axes. 

\begin{figure}
    \centering
    \includegraphics[width=0.45\textwidth]{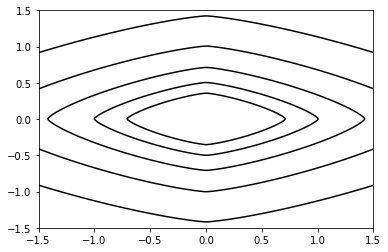}
    \includegraphics[width=0.45\textwidth]{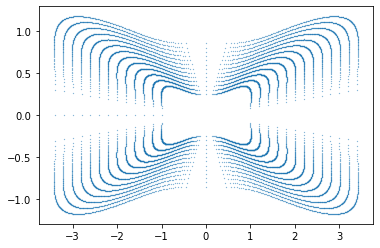}
    \caption{Contours of equal energy for the model  \eqref{eq:ex1b} with $p=1.3$ (left), and locus of the $(h_1,h_2)$-curves (right) for constant induction $|(b_1,b_2)|$.}
    \label{fig:energy}
\end{figure}

\paragraph{Energy setting.}
In the last example, we consider the corresponding interpolation in the energy framework \eqref{eq:4}. 
By elementary computations, we here obtain $\hat b_1(w) = \sqrt{2w}/c_1$ and $\hat b_2(w)=\sqrt{2w}/c_2$, and from the calculataions in Lemma~\ref{lem:5},  we further see hat 
\begin{align} \label{eq:ex1b}
w(b_1,b_2)=\tfrac{1}{2}\|(c_1 b_1, c_2 b_2\|_p^2.
\end{align}
In Figure~\ref{fig:energy}, we depict the resulting contours of equal energy and the locus of the $(h_1,h_2)$-curves for the exponent $p=13/10$. 
The illustrationss indicate, that the material bahavior is exactly the same as that of the corresponding coenergy model depicted in Figure~\ref{fig:coenergy}, which is not by coincidence:
In fact the function $w_*(h_1,h_2)$ in \eqref{eq:ex1} with exponent $n$ here is exactly the convex-conjugate of the function $w(b_1,b_2)$ with conjugate exponent $p=n/(n-1)$, and consequently $b=\nabla_h w_*(h)$ is equivalent to $h=\nabla_b w(h)$; compare with \eqref{eq:1}. 
This observation should in fact be valid for the  special setting considered in Lemma~\ref{lem:5}, but can most probably again not be transferred to the more general case.

\begin{figure}
    \centering
    \includegraphics[width=0.45\textwidth]{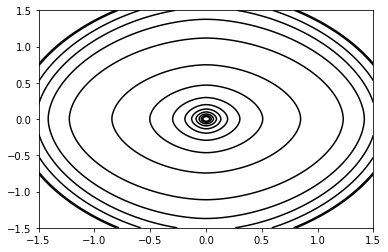}
    \includegraphics[width=0.45\textwidth]{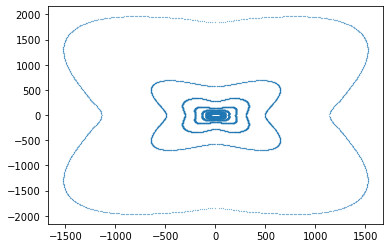}
    \caption{Contours of equal energy (left) and locus of $(h_1,h_2)$ curves for constant induction $|(b_1,b_2)|$ obtained by the implicit energy interpolation model \eqref{eq:4}.}
    \label{fig:real}
\end{figure}

\paragraph*{Example~2.}
For the last test case, we utilize the data taken from \cite[Table~1]{martin16} and use the implicit energy model \eqref{eq:4} to construct the corresponding energy functional $w(b_1,b_2)$. 
In Figure~\ref{fig:real}, we again display the contours of equal energy and the locus of the $(h_1,h_2)$-curves for constant strength of induction $|(b_1,b_2)|$.
Note that the exponent $p(w)$ in the data is not constant but close to $2$, and hence the anisotropy is not very strongly pronounced here. Nevertheless, one can again clearly see the appearance of a hard direction of magnetization different from the principal axes; also compare with the results in \cite{martin16}.

\section{Concluding discussion}

The coenergy model of Pera et al. \cite{pera93,pera94} uses measurements of magnetic material response in rolling and transverse direction in order to implicitly define a coenergy functionel $w_*(h_1,h_2)$ via solution of a scalar equation~\eqref{eq:3}. 
The underlying implicit interpolation procedure has some favorable properties, namely, it \\[1ex]
(i) exactly reproduces measured material behavior along the principal axes;\\[1ex]
(ii) implicitly defines an underlying coenergy functional $w_*(h_1,h_2)$ with convex contours of equal coenergy. \\[1ex]
(iii) allows to predicts a direction of hard magnetization different from the principal axes. \\[1ex]
By a similar relation \eqref{eq:4}, a corresponding energy functional $w(b_1,b_2)$ has been defined which shares the main feature of the coenergy model.
Based on the energy or coenergy functionals, one can use elementary relations \eqref{eq:1}--\eqref{eq:2} to compute anisotropic nonlinear material behavior for arbitrary excitations.

\paragraph{Main findings.}
We investigated mathematical properties of the implicitly defined material models, namely well-posedness as well as smoothness and convexity of the resulting (co)energy functionals. The main results of our theoretical considerations can be summarized as follows: \\[1ex]
(i) Well-posedness of the model, i.e., unique solvability of the underlying equations \eqref{eq:3} resp. \eqref{eq:4}, can be guaranteed for constant positive exponents. For variable exponents $n(w_*)$ or $p(w)$, this property may be lost. \\[1ex]
(ii) For appropriate exponents, smoothness of the (co)energy functionals can be ensured and \eqref{eq:1}--\eqref{eq:2} can be used to formally define corresponding material laws. 
\\[1ex]
(iii) Explicit formulas can be derived for computing the material relations required for finite element analysis, without the need for numerical differentiation. \\[1ex]
(iv) Convexity of the (co)energy functionals can be guaranteed in special situations, but may be lost in general.%

\paragraph{Conclusions.}
While the implicit material models \eqref{eq:3} and \eqref{eq:4} have some favorable properties and, in principle, can be evaluated efficiently, their well-posedness and convexity can only be guaranteed in special situations. 
This may be a main obstacle for the direct application of these implicit material models in the formulation of problems in nonlinear magnetostatics and their finite element analysis. 
The explicit form \eqref{eq:special}  however provides a suitable ansatz for the coenergy functional, which can also be extended easily to the coenergy framework, both extending the usual isotropic nonlinear models. 
Alternative explicit respresentations of (co)energy functionals, e.g., by by splines, are also of interest. 
The implicit material models investigated in this paper may then be valuable as priors in the fitting of such general respresentations to measurements via nonlinear least squares. Convexity of the resulting (co)energy functionals should now be guaranteed by constraints in the fitting process. The detailed investigation of such approaches is left as a topic for future research. 

{\small 
\bigskip 

\noindent 
\textbf{Acknowledgements.}
The first author is supported by FWF via the grant SFB~F9002 in the framework for the CRC "Computational Electric Machine Laboratory".

}
\bibliographystyle{spmpsci}
\bibliography{anisotropic}

\end{document}